\documentclass[11pt, a4paper]{article}
\usepackage{times}
\usepackage{a4wide}
\usepackage[dvips, hyperindex]{hyperref}
\usepackage[british]{babel}
\usepackage{enumerate}
\usepackage{amsmath, amscd, amsfonts, amssymb, latexsym, theorem}
\usepackage[T1]{fontenc}
\usepackage[latin1]{inputenc}

\theoremstyle{change}
{\theorembodyfont{\itshape}   \newtheorem{thm}{Theorem.}[section]}
{\theorembodyfont{\itshape}   \newtheorem{lem}[thm]{Lemma.}}
{\theorembodyfont{\itshape}   }
{\theorembodyfont{\itshape}   \newtheorem{rem}[thm]{Remark.}}
{\theorembodyfont{\itshape}   \newtheorem{prop}[thm]{Proposition.}}
{\theorembodyfont{\itshape}   \newtheorem{cor}[thm]{Corollary.}}

\newcommand{\QQ}{\mathbb{Q}}

\newcommand{\ZZ}{\mathbb{Z}}
\newcommand{\FF}{\mathbb{F}}
\newcommand{\NN}{\mathbb{N}}
\newcommand{\Qbar}{{\overline{\QQ}}}
\newcommand{\Zbar}{{\overline{\ZZ}}}
\newcommand{\Fbar}{{\overline{\FF}}}
\newcommand{\rhobar}{{\overline{\rho}}}
\newcommand{\GL}{\mathrm{GL}}
\newcommand{\SL}{\mathrm{SL}}
\newcommand{\PGL}{\mathrm{PGL}}
\newcommand{\PSL}{\mathrm{PSL}}
\newcommand{\PSp}{\mathrm{PSp}}
\newcommand{\Frob}{\mathrm{Frob}}
\newcommand{\Gal}{\mathrm{Gal}}
\newcommand{\Ind}{\mathrm{Ind}}
\newcommand{\Tr}{\mathrm{tr}}
\newcommand{\WD}{\mathrm{WD}}
\newcommand{\sems}{\mathrm{ss}}
\newcommand{\proj}{\mathrm{proj}}
\newcommand{\pf}{{\bf Proof. }}
\newcommand{\qed}{\hspace* {.5cm} \hfill $\Box$}
\newcommand{\mat}[4]{
 \left(  \begin{smallmatrix} #1 & #2 \\ #3 & #4 \end{smallmatrix} \right)}

\begin{document}

\title{On projective linear groups over finite fields
as Galois groups over the rational numbers}
\author{Gabor Wiese}
\maketitle
 
\begin{abstract}
  Ideas from Khare's and Wintenberger's article on the
  proof of Serre's conjecture for odd conductors are used to establish
  that for a fixed prime~$l$ infinitely many of the groups
  $\PSL_2(\FF_{l^r})$ (for $r$ running) occur as Galois groups over
  the rationals such that the corresponding number fields are
  unramified outside a set consisting of $l$, the infinite place and
  only one other prime.
\end{abstract}

\section{Introduction}

The aim of this article is to prove the following theorem.

\begin{thm}\label{mainthm}
  Let $l$ be a prime and $s$ a positive number. Then there exists a
  set~$T$ of rational primes of positive density such that for each $q
  \in T$ there exists a modular Galois representation
  $$ \overline{\rho}: G_\QQ \to \GL_2(\Fbar_l)$$
  which is unramified outside $\{\infty,l,q\}$ and whose projective
  image is isomorphic to $\PSL_2(\FF_{l^r})$ for some $r>s$.
\end{thm}

\begin{cor}\label{corone}
  Let $l$ be a prime. Then for infinitely many positive integers~$r$
  the groups $\PSL_2(\FF_{l^r})$ occur as a Galois group over the
  rationals.  \qed
\end{cor}

Using $\SL_2(\FF_{2^r}) \cong \PSL_2(\FF_{2^r})$ one obtains the
following reformulation for $l=2$.

\begin{cor}\label{cortwo}
  For infinitely many positive integers~$r$ the group
  $\SL_2(\FF_{2^r})$ occurs as a Galois group over the rationals.
  \qed
\end{cor}

This contrasts with work by Dieulefait, Reverter and Vila
(\cite{D1},\cite{D2}, \cite{RV} and~\cite{DV}) who proved that the
groups $\PSL_2(\FF_{l^r})$ and $\PGL_2(\FF_{l^r})$ occur as Galois
groups over~$\QQ$ for fixed (small)~$r$ and infinitely many
primes~$l$.

Dieulefait (\cite{D3}) has recently obtained a different, rather
elementary proof of Corollary~\ref{corone} under the assumption $l\ge
5$ with a different ramification behaviour, namely $\{\infty,2,3,l\}$.
His proof has the virtue of working with a family of modular forms
that does not depend on~$l$.

In the author's PhD thesis~\cite{Thesis} some computational evidence
on the statement of Corollary~\ref{cortwo} was exhibited. More
precisely, it was shown that all groups $\SL_2(\FF_{2^r})$ occur as
Galois groups over~$\QQ$ for $1 \le r \le 77$, extending results by
Mestre (see \cite{S}, p.~53), by computing Hecke eigenforms of
weight~$2$ for prime level over finite fields of characteristic~$2$.
However, at that time all attempts to prove the corollary failed,
since the author could not rule out theoretically that all Galois
representations attached to modular forms with image contained in
$\SL_2(\FF_{2^r})$ for $r$ bigger than some fixed bound and not in any
$\SL_2(\FF_{2^a})$ for $a \mid r$, $a \neq r$, have a dihedral image.

The present paper has undergone some developments since the first
writing. However, the core of the proof has remained the same.  It
uses a procedure borrowed from the ground breaking paper~\cite{KW1} by
Khare and Wintenberger. The representations that we will construct in
the proof are almost - in the terminology of~\cite{KW1} -
good-dihedral. The way we make these representation by level raising
is adopted from a part of the proof of~\cite{KW1}, Theorem~3.4.
Although the paper~\cite{KW1} is not cited in the proof, this paper
owes its mere existence to it.

Meanwhile, the results of the present article have been generalised to
the groups $\PSp_{2n}(\FF_{l^r})$ for arbitrary~$n$ by Khare, Larsen
and Savin (see~\cite{KLS}).  Their paper~\cite{KLS} also inspired a
slight strengthening of the main result of the present article
compared to an early version.

\subsection*{Acknowledgements}

The author would like to thank Bas Edixhoven for very useful
discussions and Alexander Schmidt for an elegant argument used in the
proof of Lemma~\ref{order}.

\subsection*{Notations}

We shall use the following notations. By $S_k(\Gamma_1(N))$ we mean
the complex vector space of holomorphic cusp forms on~$\Gamma_1(N)$
of weight~$k$. The notation $S_k(N,\chi)$ is used for the vector
space of holomorphic cusp forms of level~$N$, weight~$k$ for the
Dirichlet character~$\chi$. If $\chi$ is trivial, we write $S_k(N)$
for short.
We fix an algebraic closure~$\Qbar$ of~$\QQ$ and
for every prime~$p$ an algebraic closure~$\Qbar_p$ of~$\QQ_p$ and
we choose once and for all an embedding $\Qbar \hookrightarrow \Qbar_p$ 
and a ring surjection $\Zbar_p \to \Fbar_p$, subject to
which the following constructions are made. By $G_\QQ$ we denote the
absolute Galois group of the rational numbers. For a rational
prime~$q$, we let $D_q$ and $I_q$ be the corresponding decomposition
and inertia group of the prime~$q$, respectively.
Given an eigenform~$f \in S_k(\Gamma_1(N))$ one can attach to it by
work of Shimura and Deligne a Galois representation $\rho_{f,p}: G_\QQ
\to \GL_2(\Qbar_p)$ with well-known properties.  Choosing a lattice,
reducing and semi-simplifying, one also obtains a representation
$\rhobar_{f,p}: G_\QQ \to \GL_2(\Fbar_p)$. We denote the composition
$G_\QQ \to \GL_2(\Fbar_p) \twoheadrightarrow \PGL_2(\Fbar_p)$
by~$\rhobar_{f,p}^\proj$. All Galois representations in this paper
are continuous. By $\zeta_{p^r}$ we always mean a primitive $p^r$-th
root of unity.

\section{On Galois representations}

The basic idea of the proof is to obtain the groups
$\PSL_2(\FF_{l^r})$ as the image of some~$\rhobar_{g,l}^\proj$. In
order to determine the possible images of~$\rhobar_{g,l}^\proj$, we
quote the following well-known group theoretic result due to Dickson
(see \cite{Hu}, II.8.27).

\begin{prop}[Dickson]\label{dickson}
  Let $l$ be a prime and $H$ a finite subgroup of $\PGL_2(\Fbar_l)$.
  Then a conjugate of $H$ is isomorphic to one of the following
  groups:
\begin{itemize}
\item finite subgroups of the upper triangular matrices,
\item $\PSL_2(\FF_{l^r})$ or $\PGL_2(\FF_{l^r})$ for $r \in \NN$,
\item dihedral groups $D_r$ for $r \in \NN$ not divisible by $l$,
\item $A_4$, $A_5$ or $S_4$.
\end{itemize}
\end{prop}

We next quote a result by Ribet showing that the images
of~$\rhobar_{g,l}^\proj$ for a non-CM eigenform~$g$ are ``almost
always'' $\PSL_2(\FF_{l^r})$ or $\PGL_2(\FF_{l^r})$ for some $r \in
\NN$.

\begin{prop}[Ribet]\label{Ribet}
  Let $f = \sum_{n \ge 1} a_n q^n \in S_2(N,\chi)$ 
  be a normalised eigenform of level~$N$ and some
  character~$\chi$ which is not a CM-form.

  Then for almost all primes~$p$, i.e.\ all but finitely many,
  the image of the representation
  $$ \rhobar_{f,p}^\proj: G_\QQ \to \PGL_2(\Fbar_p)$$
  attached to~$f$ is equal to either $\PGL_2(\FF_{p^r})$ or
  to $\PSL_2(\FF_{p^r})$ for some $r \in \NN$.
\end{prop}

\pf
Reducing modulo~$p$, Theorem~3.1 of~\cite{Ri} gives, for almost all~$p$,
that the image of~$\rhobar_{f,p}$ contains $\SL_2(\FF_p)$. This
already proves that the projective image is of the form
$\PGL_2(\FF_{p^r})$ or $\PSL_2(\FF_{p^r})$ by
Proposition~\ref{dickson}.
\qed
\medskip

In view of Ribet's result, there are two tasks to be solved for
proving Theorem~\ref{mainthm}. The first task is to avoid the
``exceptional primes'', i.e.\ those for which the image
is not as in the proposition. The second one is to obtain at the same
time that the field $\FF_{p^r}$ is ``big''.

We will now use this result by Ribet in order to establish the simple
fact that for a given prime~$l$ there exists a modular form of level a
power of~$l$ having only finitely many ``exceptional primes''. The
following lemma can be easily verified using e.g.\ William Stein's
modular symbols package which is part of {\sc Magma} (\cite{Magma}).

\begin{lem}\label{lemmag}
  In any of the following spaces there exists a newform without CM:
  $S_2(2^7)$, $S_2(3^4)$, $S_2(5^3)$, $S_2(7^3)$, $S_2(13^2)$.
  \qed
\end{lem}

\begin{prop}\label{propnonsolv}
Let $l$ be a prime. Put
$$N := \begin{cases}
2^7, & \textrm{if } l=2,\\
3^4, & \textrm{if } l=3,\\
5^3, & \textrm{if } l=5,\\
7^3, & \textrm{if } l=7,\\
13^2, & \textrm{if } l=13,\\
l, & \textrm{otherwise.}
\end{cases}$$
Then there exists an eigenform $f \in S_2(N)$ such that for almost
all primes~$p$, i.e.\ for all but finitely many,
the image of the attached Galois representation~$\rhobar_{f,p}$
is of the form $\PGL_2(\FF_{p^r})$ or $\PSL_2(\FF_{p^r})$ for some $r \in \NN$.
\end{prop}

\pf
If $l \in \{2,3,5,7,13\}$, we appeal to Lemma~\ref{lemmag} to get
an eigenform~$f$ without CM. If $l$ is not in that list, then
there is an eigenform in $S_2(l)$ and it is well-known that it does
not have CM, since the level is square-free.
Hence, Proposition~\ref{Ribet} gives the claim.
\qed
\medskip

We will be able to find an eigenform with image as in the
preceding proposition which is ``big enough'' by applying the
following level raising result by Diamond and Taylor. It is a
special case of Theorem~A of~\cite{DT}.

\begin{thm}[Diamond, Taylor]\label{thmdt}
Let $N \in \NN$ and let $p > 3$ be a prime not dividing~$N$.
Let $f \in S_2(N,\chi)$ be a newform such that $\rhobar_{f,p}$
is irreducible.
Let, furthermore, $q \nmid N$ be a prime such that $q \equiv -1 \mod p$
and $\Tr(\rhobar_{f,p}(\Frob_q))=0$.

Then there exists a newform $g \in S_2(Nq^2,\tilde{\chi})$ such that
$\rhobar_{g,p} \cong \rhobar_{f,p}$. \qed
\end{thm}

\begin{cor}\label{cordt}
Assume the setting of Theorem~\ref{thmdt}.
Then the following statements hold.
\begin{enumerate}[(a)]
\item The mod~$p$ reductions of $\chi$ and $\tilde{\chi}$ are equal.

\item The restriction of $\rho_{g,p}$ to~$I_q$, the inertia group
  at~$q$, is of the form $\mat {\psi}*0{\psi^q}$ with a character $I_q
  \xrightarrow{\psi} \ZZ[\zeta_{p^r}]^\times \hookrightarrow
  \ZZ_p[\zeta_{p^r}]^\times$ of order~$p^r$ for some~$r > 0$.

\item For all primes $l \neq p,q$, the restriction of $\rhobar_{g,l}$
  to~$D_q$, the decomposition group at~$q$, is irreducible and the
  restriction of $\rhobar_{g,l}$ to~$I_q$ is of the form $\mat
  {\psi}*0{\psi^q}$ with the character $I_q \overset{\psi}{\to}
  \ZZ[\zeta_{p^r}]^\times$~\hspace*{-4pt}~
  $\twoheadrightarrow \FF_l(\zeta_{p^r})^\times$ of the same
  order~$p^r$ as in~(b).

\end{enumerate}
\end{cor}

\pf
(a) This follows from $\rhobar_{g,p} \cong \rhobar_{f,p}$.

(b) As $q^2$ precisely divides the conductor of~$\rho_{g,p}$, the
ramification at~$q$ is tame and the restriction to~$I_q$ is of the
form $\mat {\psi_1}*0{\psi_2}$ for non-trivial characters $\psi_i: I_q
\to \Zbar_p^\times$. Their order must be a power of~$p$, since their
reductions mod~$p$ vanish by assumption. Due to $q \equiv -1 \mod p$,
local class field theory tells us that $\psi_i$ cannot extend
to~$G_{\QQ_q}$, as their orders would divide~$q-1$. Hence, the image
of $\rho_{g,p}|_{D_q}: D_q \to \GL_2(\Qbar_p) \to \PGL_2(\Qbar_p)$ is
a finite dihedral group.  Consequently, $\rho_{g,p}|_{D_q}$ is an
unramified twist of the induced representation $\Ind_{D_q}^{G_K}
(\psi)$ with $K$ the unramified extension of degree~$2$ of~$\QQ_q$ and
$\psi: G_K \to
\ZZ[\zeta_{p^r}]$ a totally ramified character of order~$p^r$ for
some~$r>0$. Conjugation by $\rho_{g,p}(\Frob_q)$ exchanges $\psi_1$
and~$\psi_2$. It is well-known that this conjugation also raises
to the $q$-th power. Thus, we find $\psi_2 = \psi_1^q$ and without
loss of generality $\psi = \psi_1$.

(c) With the notations and the normalisation used in~\cite{CDT}, p.~536,
the local Langlands correspondence reads (for $l \neq q$)
$$ (\rho_{g,l}|W_q)^\sems \cong \Qbar_l \otimes_{\Qbar} \WD(\pi_{g,q}),$$
if $g$ corresponds to the automorphic representation $\pi_g =
\otimes \pi_{g,p}$.  In particular, knowing by~(b) that
$\rho_{g,p}|_{I_q}$ is $\mat {\psi}*0{\psi^q}$, it follows that
$\rho_{g,l}|_{I_q}$ is also of that form. As $\psi$ is of $p$-power
order, it cannot vanish under reduction mod~$l$ if $l\neq p$. Hence,
$\rhobar_{g,l}|_{I_q}$ is of the claimed form. The irreducibility
follows as in~(b).
\qed

\section{Proof and remarks}

In this section we prove Theorem~\ref{mainthm} and comment on possible
generalisations.

\begin{lem}\label{order}
  Let $l$ be a prime and $s\in \NN$. Then there is a set of primes~$p$
  of positive density such that $p$ is $1$ mod~$4$ and such that if
  the group $\GL_2(\FF_{l^t})$ possesses an element of order~$p$ for
  some~$t \in \NN$, then $2 \mid t$ and~$t > s$.
\end{lem}

\pf
We take the set of primes~$p$ such that $p$ is split in~$\QQ(i)$ and
inert in $\QQ(\sqrt{l})$. By Chebotarev's density theorem this set has
a positive density. The first condition imposed on~$p$ means that $4$
divides the order of~$\FF_p^\times$ and the second one that the order
of the $2$-Sylow subgroup of~$\FF_p^\times$ divides the order of~$l$
in~$\FF_p^\times$. If the order of $g \in \GL_2(\FF_{l^t})$ is~$p$,
then $\FF_{l^{2t}}^\times$ contains an element of order~$p$ (an eigenvalue
of~$g$). Thus, $l^{2t}-1$ is divisible by~$p$, whence the order of~$l$
in~$\FF_p^\times$ divides~$2t$ and consequently $2$ divides~$t$. The
condition $t>s$ can be met by excluding finitely many~$p$.
\qed
\medskip

{\bf Proof of Theorem~\ref{mainthm}.}  
Let us fix the prime~$l$ and the number~$s$ from the statement of
Theorem~\ref{mainthm}. We will now exhibit a modular form~$g$ 
such that $\rhobar_{g,l}^\proj$ has
image equal to~$\PSL_2(\FF_{l^t})$ with $t>s$.
We start with the eigenform $f \in S_2(l^*)$ provided by
Proposition~\ref{propnonsolv}.
Next, we let~$p$ be any of the infinitely many primes
from Lemma~\ref{order} such that the image of $\rhobar_{f,p}^\proj$ 
is $\PSL_2(\FF_{p^r})$ or $\PGL_2(\FF_{p^r})$ with some~$r$.
We want to obtain~$g$ by level raising. The next lemma
will yield a set of primes of positive density at which the level 
can be raised in a way suitable to us.

\begin{lem}\label{lemraise}
Under the above notations and assumptions, the set of primes~$q$ such that
\begin{enumerate}[(i)]
\itemsep=0cm plus 0pt minus 0pt
\item $q \equiv -1 \mod p$,
\item $q$ splits in $\QQ(i,\sqrt{l})$ and
\item $\rhobar_{f,p}^\proj(\Frob_q)$ lies in the same conjugacy class
as $\rhobar_{f,p}^\proj(c)$, where $c$ is any complex conjugation,
\end{enumerate}
has a positive density.
\end{lem}

\pf
The proof is adapted from \cite{KW1}, Lemma~8.2.
Let $L := \QQ(\zeta_p,i,\sqrt{l})$ and $K/\QQ$ such that 
$G_K = \ker(\rhobar_{f,p}^\proj)$. Conditions (i) and~(ii) must
be imposed on the field~$L$ and Condition~(iii) on~$K$.
We know that $\Gal(K/\QQ)$ is either $\PSL_2(\FF_{p^r})$
or $\PGL_2(\FF_{p^r})$. In the former case
the lemma follows directly from Chebotarev's density theorem,
as the intersection $L \cap K$ is~$\QQ$, since $\PSL_2(\FF_{p^r})$
is a simple group. In the latter case the intersection $L \cap K = M$
is an extension of~$\QQ$ of degree~$2$.
As $p \equiv 1 \mod 4$ by assumption, $\rhobar_{f,p}^\proj(c)$
is in $\Gal(L/M) \cong \PSL_2(\FF_{p^r})$, since $\det(\rhobar_{f,p}(c))=-1$
is a square mod~$p$.
Consequently, any~$q$ satisfying Condition~(iii) is split
in $M/\QQ$. Again as $p \equiv 1 \mod 4$, complex conjugation fixes
the quadratic subfield of~$\QQ(\zeta_p)$, whence any prime~$q$ satisfying
Conditions (i) and~(ii) is also split in $M/\QQ$. Hence,
we may again appeal to Chebotarev's density theorem, proving
the lemma.
\qed
\medskip

To continue the proof of Theorem~\ref{mainthm} we let $T$ be the set
of primes provided by Lemma~\ref{lemraise}. Let $q \in T$.
Condition~(iii) assures that $\rhobar_{f,p}(\Frob_q)$ has trace~$0$,
since it is a scalar multiple of the matrix representing complex
conjugation, i.e.~$\mat 100{-1}$. Hence, we are in the situation to
apply Theorem~\ref{thmdt} and Corollary~\ref{cordt}. We, thus, get an
eigenform $g \in S_2(l^*q^2,\chi)$ with $\chi$ a Dirichlet character
of order a power of~$p$ (its reduction mod~$p$ is trivial) such that
$\rhobar_{g,l}|_{D_q}$ is irreducible and $\rhobar_{g,l}|_{I_q}$ is of
the form $\mat {\psi}*0{\psi^q}$ with $\psi$ a non-trivial character
of order~$p^r$.  Hence, in particular, the image 
$\rhobar_{g,l}^\proj(G_\QQ)$ contains an element of order~$p$, as
$\mat {\psi}*0{\psi^q}$ cannot be scalar due to $p \nmid q-1$.

We next show that $\rhobar_{g,l}^\proj(G_\QQ)$ is not a dihedral
group.  It cannot be cyclic either because of irreducibility. If it
were dihedral, then $\rhobar_{g,l}^\proj$ would be a representation
induced from a character of a quadratic extension~$R/\QQ$. If $R$ were
ramified at~$q$, then $\rhobar_{g,p}^\proj(I_q)$ would have even
order, but it has order a power of~$p$. So, $R$ would be either
$\QQ(\sqrt{l})$ or $\QQ(\sqrt{-l})$. As by Condition~(ii) $q$~would
split in~$R$, this implies that $\rhobar_{g,p}^\proj|_{D_q}$ would
be reducible, but it is irreducible. Consequently,
$\rhobar_{g,l}^\proj(G_\QQ)$ is either $\PSL_2(\FF_{l^t})$ or
$\PGL_2(\FF_{l^t})$ for some~$t$. By what we have seen above,
$\PGL_2(\FF_{l^t})$ then contains an element of order~$p$. Hence, so
does $\GL_2(\FF_{l^t})$ and the assumptions on~$p$ imply $t>s$ and $2
\mid t$.

We know, moreover, that the determinant of $\rhobar_{g,l}(\Frob_w)$ for any prime
$w \nmid lq$ is $w^{k-1} \chi(w)$ for some fixed~$k$ (the
Serre weight of $\rhobar_{g,l}$). As $\chi(w)$ is
of $p$-power order, $\FF_{l^t}$ contains a square root of it. The
square roots of elements of~$\FF_l^\times$ are all contained
in~$\FF_{l^2}$ and thus also in $\FF_{l^t}$, as $t$ is even.
Hence, $\rhobar_{g,l}^\proj(\Frob_w) \in \PSL_2(\FF_{l^t})$. As every
conjugacy class contains a Frobenius element, the proof of Theorem~\ref{mainthm}
is finished.
\qed

\begin{rem}
  One can develop the basic idea used here further, in particular, in
  order to try to establish an analogue of Theorem~\ref{mainthm} such
  that the representations ramify at a given finite set of primes~$S$
  and are unramified outside $S \cup \{\infty,l,q\}$.
\end{rem}

\begin{rem}
  It is desirable to remove the ramification at~$l$. For that, one
  would need that $\rhobar_{g,l}^\proj$ is unramified at~$l$. This,
  however, seems difficult to establish.
\end{rem}

\vspace*{.2cm}
\noindent Gabor Wiese\\
Institut f�r Experimentelle Mathematik, Ellernstra�e 29, 45326 Essen, Germany\\
E-mail: {\tt gabor@pratum.net}, Web page: {\tt http://maths.pratum.net/}

\end{document}